\documentclass{amsart}

\usepackage{amsmath,amssymb,amscd,amsfonts}
\usepackage[all]{xy}

\newtheorem{theorem}{Theorem}
\newcommand{\bt}{\begin{theorem}}
\newcommand{\et}{\end{theorem}}
\newtheorem{lemma}{Lemma}
\newcommand{\bl}{\begin{lemma}}
\newcommand{\el}{\end{lemma}}
\newtheorem{corollary}{Corollary}
\newcommand{\bc}{\begin{corollary}}
\newcommand{\ec}{\end{corollary}}
\newcommand{\beq}{\begin{equation}}
\newcommand{\eeq}{\end{equation}}
\newcommand{\benum}{\begin{enumerate}}
\newcommand{\eenum}{\end{enumerate}}

\newcommand{\Z}{\ensuremath{\mathbf Z}}

\newcommand{\bmat}{\left(\begin{matrix}}
\newcommand{\emat}{\end{matrix}\right)}

\title{Generalizations  of Menon's arithmetic identity}
\author{Melvyn B. Nathanson}
\address{Department of Mathematics\\Lehman College (CUNY)\\Bronx, NY 10468} 
\email{melvyn.nathanson@lehman.cuny.edu} 

\dedicatory{In memoriam Eduard Wirsing} 



\thanks{Supported in part by a grant from the PSC-CUNY Research Award Program.}

\date{\today}
\begin{document}

\begin{abstract}
Menon's identity is $\sum_{a \in A} (a-1,m) = d(m) \varphi(m)$, where $A$ is a reduced set of residues modulo $m$.  This paper contains elementary proofs of  generalizations of this result. 
\end{abstract} 
\maketitle

\section{Menon's identity}

The divisor function $d(m)$ counts the number of positive divisors of the positive integer $m$.  
Let $(a,b)$ denote the greatest common divisor of integers $a$ and $b$ 
that are not both 0.  
For every  integer $s$, the \emph{$s$-free divisor function} $d_s(m)$ 
counts the number of positive divisors $d$ of $m$ such that $(d,s) = 1$. 
We have $d_0(m) = 1$ and $d_s(m) = d(m)$ if $(s,m) = 1$.  
In particular, $d_1(m) = d(m)$. 
Note that $d_s(m) = d_{-s}(m)$. 

For example, $d_1(12) = 6$, $d_2(12) = 2$, and $d_3(12) = 3$.  
If $p^v$ is a prime power, then 
\[
d_s(p^v) = \begin{cases}
1 & \text{if $p$ divides $s$} \\
v+1 &  \text{if $p$ does not divide $s$.}
 \end{cases}
\]

For every positive integer $m$, the Euler totient function $\varphi(m)$ counts the number of integers $a$ 
such that $1 \leq a \leq m$ and $(a,m) = 1$.  Thus,   
\[
\sum_{\substack{ a = 1 \\ (a,m) = 1}}^m 1 =   \varphi(m).
\]
Equivalently, 
\[
\sum_{\substack{ a = 1 \\ (a,m) = 1}}^m (a,m) =  d_0(m) \varphi(m).
\]
In 1965 P. Kesava Menon~\cite{meno65} proved the following beautiful arithmetic identity.

\bt                    \label{Menon:theorem:Menon}
For all positive integers $m$, 
\[
\sum_{\substack{ a = 1 \\ (a,m) = 1}}^m (a-1,m) = d_1(m) \varphi(m).
\]
\et

For example, if $m = 12$, then  
\[
\sum_{\substack{ a = 1 \\ (a,12) = 1}}^{12} (a-1,12) = 
(0,12) + (4,12) + (6,12) + (10,12) = 24
\] 
and 
\[
d_1(12)\varphi(12) = 6 \cdot 4 = 24.
\]
If $m = p$ is prime, then 
\[
\sum_{\substack{ a = 1 \\ (a,p) = 1}}^p (a-1,p) = (0,p) + \sum_{ b = 1}^{p-2} (b,p) 
= p + (p-2) = 2p-2 
\]
and 
\[
d_1(p) \varphi(p) = 2 (p-1) = 2p-2.
\]

Theorem~\ref{Menon:theorem:Menon} has applications to group theory.  
For example, I. M. Richards~\cite{rich84} used the identity to 
prove that if  $N$ is the number of cyclic subgroups in a finite group $G$ of order $m$, 
then $N \geq d(m)$ and $N = d(m)$ if and only if $G$ is cyclic. 

Menon~\cite{meno65}  gave three proofs of his identity.  The excellent survey paper 
of L\' aszl\' o T\' oth~\cite{toth21} 
contains 10 proofs, many generalizations, and a comprehensive bibliography 
of work on Menon's identity. 

A \emph{multiplicative function} is an arithmetic function $f(m)$ such that $f(m_1m_2) = f(m_1) f(m_2)$ 
for all pairs of relatively prime positive integers $m_1$ and $m_2$. 
The divisor functions $d(m)$ and $d_s(m)$ are multiplicative.

With the  $s$-free divisor function $d_s(m)$, there is the following 
extension of Menon's identity.

\bt           \label{Menon:theorem:MBN-1}
For all integers $s$ and all positive integers $m$, 
\[
\sum_{\substack{ a = 1 \\ (a,m) = 1}}^m (a-s,m) = d_s(m) \varphi(m) .
\]
\et

This is identity~(4.4) in T\' oth~\cite{toth21}.  See also~\cite{hauk05,hauk96}. 
Theorem~\ref{Menon:theorem:MBN-1} reduces to Theorem~\ref{Menon:theorem:Menon}
if $s=1$.

For example, if  $s=2$ and $m=12$, then 
\[
\sum_{\substack{ a = 1 \\ (a,12) = 1}}^{12} (a-2,12) = 
(-1,12) + (3,12) + (5,12) + (9,12) = 8
\]
and 
\[
d_2(12) \varphi(12) = 2\cdot 4 = 8.
\]

Let $k$ be a positive integer and let $a$ and $b$ be integers that are not both 0.
The \emph{$k$th power greatest common divisor} of $a$ and $b$, denoted $(a,b)_k$, 
is the largest $k$th power that divides both $a$ and $b$. 
The integers $a$ and $b$ are \emph{$k$th power relatively prime} if $(a,b)_k = 1$. 
For example, the integers 4 and 8 are third power relatively prime.  
The integers 8 and 27 are also third power relatively prime.

For every positive integer $m$, the 
\emph{Eckford Cohen totient function}\index{Eckford Cohen totient function}\index{totient function!Eckford Cohen} 
 $\varphi^{(k)}\left( m \right)$ (Cohen~\cite{cohe56}) counts the number of integers $a$ 
such that $1 \leq a \leq m^k$ and $ (a,m^k)_k=1$.  
Note that $ \varphi^{(1)}\left( m \right) =   \varphi \left( m \right)$ for all $m$.
Properties of the Eckford Cohen totient function are discussed 
in Section~\ref{Menon:section:Cohen}.  
For example, with $m = 2^k$ we see that 
\[
\left(a,\left(2^{k}\right)^2 \right)_2 = (a, 4^k)_2 = 1 
\] 
if and only if 4 does not divide $a$ and so 
\[
\varphi^{(2)}\left( 2^k \right)  = \sum_{\substack{a=1 \\ \left( a,4^k \right)_2 =1}}^{4^k} 1 
 =  4^k - 4^{k-1} = 3\cdot 4^{k-1}.
\]

In 1972, K. Nageswara Rao~\cite{rao72} applied the Eckford Cohen totient function 
to obtain the following generalization of Menon's identity.

\bt                          \label{Menon:theorem:Rao}
For all integers $s$ and all positive integers  $k$ and $m$,  
if $s$ and $m^k$ are $k$th power relatively prime, then     
\[
\sum_{\substack{a=1 \\ (a,m^k)_k = 1}}^{m^k} \left(a-s, m^k \right)_k = d(m)\varphi^{(k)}(m). 
\]
\et

For example, with $s=1$, $k=2$, and $m = 4$ we have 
\[
\sum_{\substack{a=1 \\ (a,16)_2 = 1}}^{16} \left(a-1, 16 \right)_2 
= \sum_{\substack{b=0 \\ b \neq 3, 7, 11 }}^{14} (b,16)_2 
= 36   
\]
and
\[
d(4)\varphi^{(2)}(4) = 3\cdot 12 = 36. 
\]

Theorem~\ref{Menon:theorem:Rao} reduces to 
Theorem~\ref{Menon:theorem:Menon}  (Menon's identity) if $k=s=1$ 
and Theorem~\ref{Menon:theorem:MBN-1} if $k=1$ and $(s,m) =1$.

To remove the $k$th power relatively prime condition from Rao's result, 
we introduce the divisor function $d_s^{(k)}(m)$, 
which is defined as follows:
If $p^v$ is a prime power, then   
\beq                                                \label{Menon:dsk}
d_s^{(k)}(p^v) = \begin{cases}
1 & \text{if $p^k$ divides $s$} \\
v+1 &  \text{if $p^k$ does not divide $s$.}
 \end{cases}
\eeq
Write $p^v \| m$ if $p^v$ divides $m$ and $p^{v+1}$ does not divide $m$.  
For every positive integer $m$, we define 
\[
d_s^{(k)}(m) = \prod_{p^v \| m} d_s^{(k)}(p^v). 
\]
The arithmetic function $d_s^{(k)}(m)$ is multiplicative, 
and $d_s^{(1)}(m) = d_s(m)$.

The following identity generalizes Rao's Theorem~\ref{Menon:theorem:Rao}. 

\bt                        \label{Menon:theorem:MBN-2}
For all integers $s$ and all positive integers  $k$ and $m$, 
\[
\sum_{\substack{a=1 \\ (a,m^k)_k = 1}}^{m^k} \left(a-s, m^k \right)_k 
=  d_s^{(k)}(m)  \varphi^{(k)}(m). 
\]
\et

For example, if $s=12$, $k=2$, and $m = 4$, then $\left(s, m^k \right)_k = \left(12, 16 \right)_2  = 4$.  
We have  
\[
\sum_{\substack{a=1 \\ (a,16)_2 = 1}}^{16} \left(a-12, 16 \right)_2 
=  \sum_{\substack{b=-11 \\ b \neq -8,-4,0} }^3 (b,16)_2 = 12 
\]
and 
\[
d_{12}^{(2)}(4) \varphi^{(2)}(4) = 1 \cdot 12 = 12.
\]

Theorem~\ref{Menon:theorem:MBN-2} reduces to 
Theorem~\ref{Menon:theorem:MBN-1} if $k=1$ 
and to Theorem~\ref{Menon:theorem:Rao} 
if the integers $s$ and $m^k$ are  $k$th power relatively prime.  
In Section~\ref{Menon:section:proof} we use multiplicative arithmetic functions to give a simple proof of Theorem~\ref{Menon:theorem:MBN-2}.

\section{The Eckford Cohen totient function $\varphi^{(k)}(m) $}   \label{Menon:section:Cohen} 
We have 
\[
\varphi^{(k)}(m) = \sum_{\substack{a=1 \\ (a,m^k)_k = 1}}^{m^k} 1. 
\]
Note that $(a,m^k)_k > 1$ if and only if there is an integer $d > 1$ such that $d$ divides $m$ 
and $d^k$ divides $a$.

\bl                   \label{Menon:lemma:CohenTotient}
The Eckford Cohen totient function 
$\varphi^{(k)}(m) $ is a multiplicative arithmetic function.  For every positive integer $m$, 
\[
 \varphi^{(k)}\left( m \right) =  m^k \prod_{p|m} \left( 1 - \frac{1}{p^k} \right).
\]
If $p$ is a prime and $v$ is a positive integer, then 
\[
 \varphi^{(k)}\left( p^v \right) =  p^{vk}- p^{(v-1)k}.
\]
\el

\begin{proof}
For every  divisor $d$ of $m$, the number of positive  integers up to $m^k$ that are divisible 
by $d^k$ is $m^k/d^k$.  In particular, if $p_{i_1} p_{i_2} \cdots p_{i_{\ell}}$ is a product 
of distinct primes that divide $m$, then the number of positive  integers up to $m^k$ 
that are divisible by $(p_{i_1} p_{i_2} \cdots p_{i_{\ell}})^k$ 
is $m^k/(p_{i_1} p_{i_2} \cdots p_{i_{\ell}})^k$.

Let $p_1,\ldots, p_r$ be the distinct primes that divide $m$.  
Let $1 \leq a \leq m^k$.  We have $(a, m^k)_k = 1$ if and only if $a$ is not divisible 
by $p_i^k$ for all $i \in \{1,2,\ldots, r\}$.  
The inclusion-exclusion principle implies that 
\begin{align*} 
 \varphi^{(k)}\left( m \right) 
 & = m^k - \sum_{i_1=1}^r \frac{m^k}{p_{i_1}^k}
 + \sum_{\substack{ i_1,i_2=1 \\ i_1 < i_2 }}^r   \frac{m^k}{p_{i_1}^kp_{i_2}^k} 
 - \sum_{\substack{ i_1,i_2,i_3=1 \\ i_1 < i_2 < i_3 }}^r   \frac{m^k}{p_{i_1}^kp_{i_2}^k p_{i_3}^k} 
 + \cdots \\
& =  m^k \prod_{i=1}^r \left( 1 - \frac{1}{p_i^k} \right).
\end{align*}
This formula implies that the function $\varphi^{(k)}(m) $ is multiplicative.  
The expression for $\varphi^{(k)}(p^v)$ follows immediately.  
This completes the proof. 
\end{proof}

\bl                    \label{Menon:lemma:kth-gcd}
If $a\equiv b \pmod{m^k}$, then $(a,m^k)_k = (b,m^k)_k$. 
If one element of a congruence class modulo $m^k$ is $k$th power relatively prime to $m^k$, 
then every element in the congruence class is $k$th power relatively prime to $m^k$. 
\el

\begin{proof}
If $a = b + qm^k$, then $d^k$ divides both $a$ and $m^k$ if and only if 
$d^k$ divides both $b$ and $m^k$.  It follows that $(a,m^k)_k = (b,m^k)_k$. 
This completes the proof. 
\end{proof}

A \emph{$k$th power reduced  set of residues modulo $m$} is a set of 
$ \varphi^{(k)}\left( m \right)$ integers, one from each congruence class modulo $m^k$ 
whose elements are $k$th power relatively prime to $m^k$. 
The standard $k$th power reduced  set of residues modulo $m$ is 
\[
 \left\{a: 1 \leq a \leq m^k \text{ and } (a,m^k)_k = 1 \right\}.
\]
Lemma~\ref{Menon:lemma:kth-gcd} implies that 
if $A$ is a $k$th power reduced  set of residues modulo $m$, then 
\[
\varphi^{(k)}(m) = \sum_{a \in A} \left( a,m^k \right)_k.
\]

\bl            \label{Menon:lemma:A1A2set} 
Let $m_1$ and $m_2$ be relatively prime positive integers.  
If $A_1$ is a $k$th power reduced  set of residues modulo $m_1$ and 
$A_2$ is a $k$th power reduced  set of residues modulo $m_2$, then the set 
\[
A = \{a_1 m_2^k + a_2 m_1^k : a_1 \in A_1 \text{ and } a_2 \in A_2\}
\]
is a $k$th power reduced  set of residues modulo $m_1m_2$.
\el

\begin{proof}
Let $a = a_1 m_2^k + a_2 m_1^k \in A$ and let $d^k = (a,m_1^km_2^k)_k$.  
Suppose that $d > 1$.  If the prime $p$  divides $d$, then $p^k$ divides $d^k$ 
and $d^k$ divides $m_1^km_2^k$, and so $p^k$ divides $m_1^km_2^k$.
If $p$ divides $m_1$, then $p$ does not divide $m_2$ and so $p^k$ divides $m_1^k$.  
Because  $p^k$ divides $a =a_1 m_2^k + a_2 m_1^k$, 
it follows that $p^k$ divides $a_1 m_2^k$ and so $p^k$ divides $a_1$.  
Thus, $1 = (a_1,m_1^k)_k \geq p^k$, 
which is absurd.  Therefore, $(a,m_1^km_2^k)_k = 1$ and every element of $A$ 
is $k$th power relatively prime to $m_1m_2$.  

Let $a_1,a'_1 \in A_1$ and $a_2,a'_2 \in A_2$.  If
\[
 a_1 m_2^k + a_2 m_1^k  \equiv  a'_1 m_2^k + a'_2 m_1^k  \pmod{m_1^k m_2^k}
\]
then  
\[
 a_1 m_2^k \equiv  a_1 m_2^k + a_2 m_1^k  \equiv  a'_1 m_2^k + a'_2 m_1^k  
 \equiv  a'_1 m_2^k  \pmod{m_1^k 
 }.
\]
Because $(m_1,m_2)=1$, we have $a_1 \equiv  a'_1   \pmod{m_1^k}$ and so $a_1 = a'_1$.   
Similarly, $a_2 =  a'_2$.   
Thus, the elements of $A$ belong to pairwise distinct congruence classes 
modulo $m_1^km_2^k$.  

By Lemma~\ref{Menon:lemma:CohenTotient}, 
the Eckford Cohen totient function is multiplicative, 
and so  
\[
\varphi^{(k)}(m_1m_2) =  \varphi^{(k)}(m_1)\varphi^{(k)}(m_2) = |A_1| \ |A_2| = |A| 
\]
and the set $A$ is a $k$th power reduced set of residues modulo $m_1m_2$.  
This completes the proof. 
\end{proof}

\section{Pillai's sum functions $P^{(k)}(m)$} 

Let $m$ and $k$ be positive integers.  
The \emph{gcd sum function}\index{gcd sum function}, 
also called \emph{Pillai's arithmetic function}\index{Pillai's arithmetic function},  is 
\[
P(m) =  \sum_{a=1}^m (a,m).
\]
The \emph{$k$th power gcd sum function}, that is, the $k$th power analogue 
of Pillai's arithmetic function, is 
\[
P^{(k)}(m) = \sum_{a=1}^{m^k} (a,m^k)_k.
\]
Note that $P^{(1)}(m) = P(m)$.  

\bl           \label{Menon:lemma:Pillai} 
For every positive integer $m$, 
\[
P^{(k)}(m) = \sum_{d|m} d^k \varphi^{(k)}\left(\frac{m}{d}\right). 
\]
If $p^v$ is a prime power, then 
\[
P^{(k)}\left( p^v \right) =   (v+1) p^{vk} - v p^{(v-1)k}.
\]
\el

\begin{proof}
Write  $a = bd^k$ if $d^k$ divides $a$. 
We have 
\begin{align*}
P^{(k)}(m) & =  \sum_{a=1}^{m^k} (a,m^k)_k 
= \sum_{d|m} d^k \sum_{\substack{a=1 \\ d^k|a \\ \left( a/d^k,  ( m/d )^k \right)_k=1}}^{m^k} 1 \\
& = \sum_{d|m} d^k \sum_{\substack{b=1 \\ (b,( m/d )^k)_k =1}}^{( m/d )^k} 1  
= \sum_{d|m} d^k \varphi^{(k)}\left(\frac{m}{d}\right). 
\end{align*}
The set of $k$th power divisors of the prime power $p^{vk}$ is 
$\left\{ p^{ik} : i = 0,1,2,\ldots, v \right\}$. 
Applying Lemma~\ref{Menon:lemma:CohenTotient}, we obtain  
\begin{align*}
P^{(k)}\left( p^v \right) 
& = \sum_{i=0}^v p^{ik} \varphi^{(k)}\left(p^{v-i} \right) 
= p^{vk} +  \sum_{i=0}^{v-1} p^{ik}  \left(p^{(v-i)k} - p^{(v-i -1)k} \right) \\
& = p^{vk} +  \sum_{i=0}^{v-1}   \left(p^{vk} - p^{(v-1)k} \right) = p^{vk} +  v  \left(p^{vk} - p^{(v-1)k} \right) \\
& =   (v+1) p^{vk} - v p^{(v-1)k}.
\end{align*}
This completes the proof.  
\end{proof}

\section{Proof of Theorem~\ref{Menon:theorem:MBN-2}}  \label{Menon:section:proof}
 
We use the following properties of the  $k$th power greatest common divisor.

\bl                   \label{Menon:lemma:gcdk}
Let $k$ and $m$ be positive integers, 
and let  $A$ be a $k$th power reduced set of residues modulo $m$.   
\benum
\item[(i)]
For all integers $a$, the arithmetic function 
$m \mapsto  (a,m^k)_k$ is multiplicative.  

\item[(ii)]
For all integers $\ell$ and $s$, if  $(\ell,m) = 1$, then 
\[
\sum_{a\in A} (a \ell - s,m^k)_k = \sum_{a\in A} (a - s,m^k)_k. 
\]
\eenum 
\el

\begin{proof}
(i)  Fix the integer $a$.  Let $m_1$ and $m_2$ be relatively prime positive integers.  
The positive integer $(a,m_1^k)_k$ divides $m_1^k$ and 
the positive integer $(a,m_2^k)_k$ divides $m_2^k$ and so 
$(a,m_1^k)_k$ and  $(a,m_2^k)_k$ are relatively prime.  
The integers $(a,m_1^k)_k$ and  $(a,m_2^k)_k$ 
also divide $(a, (m_1m_2)^k)_k$ and so 
$(a,m_1^k)_k (a,m_2^k)_k$ divides $(a, (m_1m_2)^k)_k$. 

Let $d^k = (a, (m_1m_2)^k)_k$. 
Because $(m_1,m_2)=1$ and $d^k$ divides $m_1^k m_2^k$, 
there exist unique positive integers $d_1$ and $d_2$ such that 
$d^k = d_1^k d_2^k$, where $d_1^k$ divides $m_1^k$ and $d_2^k$ divides $m_2^k$.  
It follows that $d_1^k$ divides $(a,m_1^k)_k$ and $d_2^k$ divides $(a,m_2^k)_k$, 
and so  $ (a, (m_1m_2)^k)_k = d^k = d_1^kd_2^k$ divides $(a,m_1^k)_k (a,m_2^k)_k$. 
Therefore, $(a, (m_1m_2)^k)_k = (a,m_1^k)_k (a,m_2^k)_k$. 
This proves the multiplicativity of the function 
$m \mapsto  (a,m^k)_k$. 

(ii)  Let $A$ be a $k$th power reduced set of residues modulo $m$.  Let $a \in A$ 
and let $\ell$ be an integer that is relatively prime to $m$.  
If 
\[
(a\ell,m^k)_k = d^k 
\]
then $d^k$ divides $a\ell$ and $d^k$ divides $m^k$.   
It follows from $(\ell,m) = 1$ that $(d^k,\ell) = 1$ and so $d^k$ divides $a$.  
From $(a,m^k)_k = 1$ we deduce that $d=1$ and so $(a\ell,m^k)_k = 1$ 
for all $a \in A$.  

If $a,a' \in A$ and $a \ell \equiv a' \ell \pmod{m^k}$, then 
 $a \equiv a' \pmod{m^k}$ and so $a = a'$.  Therefore,  $\{a\ell:a\in A\}$
is a $k$-th power reduced set of residues modulo $m$.  
It follows that the sets 
$\{a\ell - s:a\in A\}$ and $\{a - s:a\in A\}$ 
represent the same set of $\varphi^{(k)}(m)$ distinct congruence classes modulo $m^k$. 
This completes the proof. 
\end{proof}

For positive integers $k$ and $m$, let $A$ be a $k$-th power reduced set 
of residues modulo $m$. 
For every integer $s$, we define 
\[
 M_s^{(k)}(m) = \sum_{\substack{a \in  A\\ (a,m^k)_k = 1}} \left(a-s, m^k \right)_k. 
\]
By Lemma~\ref{Menon:lemma:kth-gcd}, 
this independent of the choice of the set $A$. 

\bl                                \label{Menon:lemma:Msk-mult}
For all integers $s$ and all positive integers $k$, 
the arithmetic function $m \mapsto M_s^{(k)}(m)$ is multiplicative.  
\el

\begin{proof}
Let $m_1$ and $m_2$ be relatively prime positive integers.  
Let $A_1$ be a $k$th power reduced set of residues modulo $m_1$ 
and let $A_2$ be a $k$th power reduced set of residues modulo $m_2$.  
By Lemma~\ref{Menon:lemma:A1A2set}, the set 
\[
A = \{a_1 m_2^k + a_2 m_1^k : a_1 \in A_1 \text{ and } a_2 \in A_2\}
\]
is a $k$th power reduced  set of residues modulo $m_1m_2$. 
Applying Lemmas~\ref{Menon:lemma:gcdk} and~\ref{Menon:lemma:kth-gcd}, 
we obtain 
\begin{align*}
M_s^{(k)}(m_1m_2) 
& = \sum_{a\in A} \left( a-s,m_1^km_2^k \right)_k \\
& = \sum_{a_1 \in A_1} \sum_{a_2 \in A_2}  \left( a_1 m_2^k + a_2 m_1^k  -s,m_1^km_2^k \right)_k \\
& = \sum_{a_1 \in A_1} \sum_{a_2 \in A_2}  \left( a_1 m_2^k + a_2 m_1^k  -s,m_1^k  \right)_k 
\left( a_1 m_2^k + a_2 m_1^k  -s, m_2^k \right)_k\\
& = \sum_{a_1 \in A_1} \sum_{a_2 \in A_2}  \left( a_1 m_2^k  -s,m_1^k  \right)_k \left(  a_2 m_1^k  -s, m_2^k \right)_k\\
& = \sum_{a_1 \in A_1}  \left( a_1 m_2^k  -s,m_1^k  \right)_k \sum_{a_2 \in A_2}  \left(  a_2 m_1^k  -s, m_2^k \right)_k\\
& =  \sum_{a_1 \in A_1}  \left( a_1    -s,m_1^k  \right)_k \sum_{a_2 \in A_2}  \left(  a_2  -s, m_2^k \right)_k\\
& = M_s^{(k)}(m_1) M_s^{(k)}(m_2).
\end{align*} 
This completes the proof.  
\end{proof}

\bl                                  \label{Menon:lemma:Msk-pv}
For every integer $s$ and every prime power $p^v$, 
\[
M_s^{(k)}(p^v) = d_s^{(k)}(p^v)\varphi^{(k)}(p^v).
\]
\el

\begin{proof}  
The set 
\[
A = \left\{ 1,2,\ldots, p^{vk}\right\} \setminus  \left\{ip^k :  i=1,2,\ldots, p^{(v-1)k} \right\}
\] 
is a $k$th power reduced set of residues modulo $p^v$ and 
\[
|A| = \varphi^{(k)}(p^v) = p^{vk} - p^{(v-1)k}.
\]  

There are two cases.  In the first case, the integer $s$ is divisible by  $p^k$.    
If \\ 
$(a-s, p^{vk})_k > 1$, then $p^k$ divides $a-s$ and so $p^k$ divides $a$ 
and $(a,p^{vk})_k > 1$. 
Therefore, $(a,p^{vk})_k = 1$ implies $(a-s,p^{vk})_k = 1$.  

We have $d_s^{(k)}(p^v) = 1$ by~\eqref{Menon:dsk}. 
By Lemma~\ref{Menon:lemma:CohenTotient}, we have  
\begin{align*}
M_s^{(k)}(p^v) 
& = \sum_{\substack{a=1 \\ (a,p^{vk})_k = 1 } }^{p^{vk}} \left( a-s, p^{vk} \right)_k 
= p^{vk} - p^{(v-1)k} \\ 
& =  \varphi^{(k)}(p^v)  = d_s^{(k)}(p^v) \varphi^{(k)}(p^v) . 
\end{align*}

 In the second case, the integer $s$ is not divisible by  $p^k$, that is, $(s,p^k)_k=1$.
If  $\left( bp^k-s, p^{vk} \right)_k > 1$ for some $b \in \Z$, then $p^k$ divides $s$, 
which is absurd.  Therefore,  
\[
\left( bp^k-s, p^{vk} \right)_k = 1
\]
for all integers $b$.    
Using the gcd sum function $P^{(k)}(m)$ from Lemma~\ref{Menon:lemma:Pillai}, we obtain  
\begin{align*}
M_s^{(k)}(p^v) 
& = \sum_{\substack{a=1 \\ (a,p^k)_k = 1 } }^{p^v} \left( a-s, p^{vk} \right)_k \\
& =  \sum_{a=1}^{p^{vk}}  \left( a-s, p^{vk} \right)_k 
- \sum_{\substack{a=1 \\ (a,p^k)_k > 1 } }^{p^v} \left( a-s, p^{vk} \right)_k \\
& =  \sum_{a=1}^{p^{vk}}  \left( a-s, p^{vk} \right)_k -  \sum_{b=1}^{p^{(v-1)k}}  \left( bp^k-s, p^{vk} \right)_k \\ 
& =  \sum_{a=1}^{p^{vk}}  \left( a, p^{vk} \right)_k -  \sum_{b=1}^{p^{(v-1)k}} 1 \\ 
& = P^{(k)}(p^v) - p^{(v-1)k} \\
& = \left(  (v+1) p^{vk}- vp^{(v-1)k} \right) - p^{(v-1)k} \\
& =  (v+1) \left( p^{vk}- p^{(v-1)k} \right)  \\
& = d_s^{(k)}  \left( p^{v} \right) \varphi^{(k)}  \left( p^{v} \right). 
\end{align*}
This completes the proof. 
\end{proof}

The arithmetic functions $M_s^{(k)}(m)$, $d_s(m)$, and $\varphi^{(k)}(m)$ 
are multiplicative.  Therefore, 
\begin{align*}
M_s^{(k)}(m) 
&  = \prod_{p^v \| m} M_s^{(k)} \left( p^v \right) 
=  \prod_{p^v \| m}  d_s^{(k)}(p^v)  \varphi^{(k)}(p^v)  \\ 
&=  d_s^{(k)}(m)  \varphi^{(k)}(m). 
\end{align*}
This completes the proof of Theorem~\ref{Menon:theorem:MBN-2}. 
\\

{\emph{Acknowledgements}. I thank L\' aszl\' o T\' oth 
and Kevin O'Bryant for helpful remarks and references.

\end{document}